\documentclass[12pt]{amsart}
\usepackage{amsfonts,amssymb,amsmath,latexsym,array,amsthm}
\usepackage{fullpage}
\usepackage{verbatim,euscript}
\usepackage{pslatex}

\input {cyracc.def}

\usepackage{color}
\definecolor{darkblue}{rgb}{0,0,.5}
\usepackage[colorlinks=true,linkcolor=blue,urlcolor=blue,citecolor=darkblue]{hyperref}

\theoremstyle{definition}
\newtheorem{df}{Definition}

\theoremstyle{plain}
\newtheorem{thm}{Theorem}
\newtheorem{lm}[df]{Lemma}

\newtheorem{prop}[df]{Proposition}

\theoremstyle{remark}
\newtheorem{ex}[df]{Example}
\newtheorem{rmk}[df]{Remark}

\newcommand{\cled }{{\noindent \bf Corollary.}}

\newcommand{\um}{{\cdot}}

\newcommand{\eus}{\EuScript}
\newcommand{\gt}{\mathfrak}
\newcommand{\cp}{\mathbb C}
\newcommand{\rl}{\mathbb R}

\newcommand{\fk}{\mathbb F}

\newcommand{\ad}{{\rm ad }}

\newcommand{\Sp}{{\rm Sp}}

\newcommand{\Lie}{{\rm Lie\,}}

\newcommand{\ind}{{\rm ind\,}}
\newcommand{\chr}{{\rm char\,}}

\newcommand {\rk}{{\mathrm{rk\,}}}
\newcommand{\cP}{{\mathcal P}}
\newcommand{\cS}{{\mathcal S}}
\newcommand{\cA}{{\mathcal A}}

\newfam\Bbbfam\newfam\eufam\newfam\eusfam%
\font\Bbbfont=msbm10 scaled 1200%
\font\Bbbsmallfont=msbm8%
\textfont\Bbbfam=\Bbbfont\scriptfont\Bbbfam=\Bbbsmallfont
\font\euzw=eufm10 scaled 1200%
\font\euac=eufm7 scaled 1200%
\font\euacc=eufm7 scaled 1000%
\textfont\eufam=\euzw\scriptfont\eufam=\euac%
\scriptscriptfont\eufam=\euacc%
\font\euszw=eusm10 scaled 1200%
\font\eusac=eusm7 scaled 1200%
\font\eusacc=eusm7 scaled 1000%
\textfont\eusfam=\euszw\scriptfont\eusfam=\eusac%
\scriptscriptfont\eusfam=\eusacc%

\begin{document}
\hfill {\scriptsize   December 31, 2017} 
\vskip1ex

\title[Complete families of commuting functions]{Complete
families of commuting functions \\ for coisotropic Hamiltonian  actions}
\author{Ernest B. Vinberg}
\address[E.Vinberg]{Moscow State University,
Faculty of Mechanics and Mathematics,
Vorob'evy gory,
\ 119992 Moscow, Russia}
\email{vinberg@zebra.ru}
\author{Oksana S. Yakimova}
\address[O.Yakimova]{Institut f\"ur Mathematik, Friedrich-Schiller-Universit\"at Jena, Jena, 07737, Deutschland}
\email{oksana.yakimova@uni-jena.de}
\thanks{The first author is partially supported by the RFBR Grant 16-01-00818;
the second author is partially supported by the Graduiertenkolleg GRK~1523 ``Quanten- und Gravitationsfelder".}
\keywords{Symplectic variety, Poisson algebra, Coisotropic action, Coadjoint representation}
\subjclass[2010]{17B63, 53D17}
\maketitle

\begin{abstract}
Let $G$ be an algebraic group defined over a field $\fk$ of characteristic zero with $\gt g=\Lie G$. The dual space $\gt g^*$ equipped with the Kirillov-Kostant bracket is a Poisson variety and each irreducible $G$-invariant subvariety $X\subset \gt g^*$ carries the induced Poisson structure. We prove that there is a set 
$\{f_1,...,f_l\}\subset \fk[X]$
of algebraically independent polynomial functions, which pairwise commute with respect to the Poisson bracket,  such that $l=(\dim X+{\rm tr.deg}\, \fk(X)^G)/2$. We also discuss several applications of this result to complete integrability of Hamiltonian systems on symplectic Hamiltonian $G$-varieties of corank zero and $2$.
\end{abstract}

\section*{Introduction}
In this paper, we study Hamiltonian actions of algebraic groups on affine varieties focusing on the non-reductive case.  The ground field $\fk$ is assumed to be of characteristic zero, but not necessarily  algebraically closed.
Let us start with main definitions in the general algebraic setting.

\begin{df} Let $\cA$ be a commutative associative $\fk$-algebra
equipped with an additional anticommutative bilinear operation
$\{\,\,,\,\,\}:\,{\mathcal A}\times{\mathcal A}\to{\mathcal
A}$ called a {\it Poisson bracket} such that
$$
\begin{array}{c}
\{a,bc\}=\{a,b\}c+b\{a,c\},\\
\{a,\{b,c\}\}+\{b,\{c,a\}\}+\{c,\{a,b\}\}=0\\
\end{array}
$$
for all $a,b,c\in{\cA}$. Then $\cA$ is called a {\it Poisson algebra}.  An ideal  $I\lhd \cA$ is said to be {\it Poisson} if $\{I,\cA\}\subset I$; a homomorphism $\varphi:\,{\cA}\to{\mathcal B}$ of Poisson algebras is said to be {\it Poisson} if
$\varphi(\{x,y\})=\{\varphi(x),\varphi(y)\}$ for all $x,y\in{\cA}$. The {\it Poisson centre} of ${\mathcal A}$ is the Poisson subalgebra
$ {\eus Z}{\mathcal A}:=\{a\in{\mathcal A} \mid \{a,{\mathcal A}\}=0\}$.
A subalgebra  ${\mathcal B}\subset {\mathcal A}$
is said to be {\it Poisson-commutative} if
$\{\mathcal B, \mathcal B\}=0$.
\end{df}

Let ${\cP}$ be a Poisson algebra. Assume  that $\cP$ has no zero-divisors and
${\rm tr.deg}\,{\cP}< \infty$.  Let ${\rm Der}(\cP, {\rm Quot}\,{\cP})$ stand for the set of all ${\rm Quot}\,{\cP}$-valued derivations of the algebra $\cP$ regarded just as a commutative associative algebra. This is a linear space over ${\rm Quot}\,{\cP}$ of dimension ${\rm tr.deg}\,{\cP}$. Each $\varphi\in\cP$ gives rise to a derivation $\ad(\varphi)$, where $\ad(\varphi){\cdot}\psi=\{\varphi,\psi\}$ for all $\psi\in\cP$. Let $V({\cP}):=\langle\ad(\varphi)\mid \varphi\in {\cP}\rangle$ be the subspace of ${\rm Der}(\cP, {\rm Quot}\,{\cP})$ spanned by the inner derivations.
Then $\dim V({\cP})$ (the dimension over ${\rm Quot}\,{\cP}$) is said to be the {\it rank} of ${\cP}$, usually denoted by $\rk{\cP}$.

If the ground field $\fk$ is algebraically closed and the algebra ${\cP}$ is finitely generated, then ${\rm Der}(\cP, {\rm Quot}\,{\cP})$ can be viewed as the space of rational vector fields on the affine algebraic variety ${\rm Spec}\,{\cP}$. The inner derivations of $\cP$ are then interpreted as the Hamiltonian vector fields.

Next, set $\omega(\ad(\varphi), \ad(\psi)):=\{\varphi,\psi\}$. Since $\ad(\varphi)=0$ for each $\varphi\in{\eus Z}{\cP}$, $\omega$ is a non-degenerate skew-symmetric bilinear form on $V({\cP})$ over ${\rm Quot}\,{\cP}$.  Hence, in particular, $\rk{\cP}$ is even. It is not difficult to see that $V({\cP})$ and
$\omega$ do not change if we pass to the localisation  of ${\cP}$ by a 
multiplicative subset of ${\eus Z}{\cP}$.

\begin{df} \label{sympl} A Poisson algebra ${\cP}$ is said to be {\it symplectic}, if
$V({\cP})= {\rm Der}(\cP, {\rm Quot}\,{\cP})$, or, in other words,
if $\rk{\cP}={\rm tr.deg}\,{\cP}$.
\end{df}

\begin{df}\label{ham} A {\it Hamiltonian action} of a (finite-dimensional) Lie algebra $\gt q$ on a symplectic algebra $\cP$ is a linear map $\rho\!: \gt q\to \cP$ such that $\rho([\xi,\eta])=\{\rho(\xi),\rho(\eta)\}$ for all $\xi,\eta\in\gt q$ and each $p\in\cP$ is contained in an $\ad(\rho(\gt q))$-invariant finite-dimensional subspace of $\cP$.
\end{df}

In what follows, we assume that $\rho$ is injective  and consider $\gt q$ as a Lie subalgebra of $\cP$. The  Poisson subalgebra $\cP(\gt q)\subset\cP$, generated by $\gt q$, is called the {\it Noether subalgebra}.

Let $\cP$ be a symplectic algebra and $\cA\subset\cP$ a Poisson subalgebra.
Let $U(\cA)\subset V(\cP)$ be the subspace spanned over ${\rm Quot}\,{\cP}$ by the derivations $\ad(\varphi)$ with $\varphi\in\cA$.

\begin{df}\label{cois} A Hamiltonian action  $\gt q\hookrightarrow\cP$ is said to be
{\it coisotropic} if the subspace $U(\cP(\gt q))$ is coisotropic with respect to $\omega$.
\end{df}

The main result of the paper is the following theorem.

\begin{thm}\label{main-new}
For any coisotropic Hamiltonian action of a Lie algebra $\gt q$ on a symplectic algebra $\cP$, the subalgebra $\cP(\gt q)$ contains a Poisson-commutative subalgebra of transcendence degree $\frac{1}{2}\rk\cP$.
\end{thm}

With a few preparations it follows from a more geometric statement. Let $\gt g=\Lie G$ be the Lie algebra of a connected algebraic (or a Lie) group $G$, and
${\mathcal S}(\gt g)$ be the symmetric algebra of $\gt g$. Then ${\mathcal S}(\gt g)$ is a Poisson algebra and the algebra $\gt g$ acts on it in the sense of Definition~3. The same holds for any quotient of ${\mathcal S}(\gt g)$ by a $G$-invariant ideal $I\lhd {\mathcal S}(\gt g)$ 
(which is automatically a Poisson ideal).

\begin{thm} \label{main2}
Let $I\lhd{\mathcal S}(\gt g)$ be a prime $G$-invariant ideal. Set
$$
l(I):={\rm tr.deg} ( {\mathcal S}(\gt g)/I ) -\frac{1}{2}\rk ({\mathcal S}(\gt g)/I).
$$
Then there are Poisson-commuting algebraically independent functions
$f_1, \ldots, f_{l(I)}\in  {\mathcal S}(\gt g)/I$.
\end{thm}

In case $I=0$, the existence of a Poisson-commutative subalgebra
$\cA\subset\cS(\gt g)$ with ${\rm tr.deg}\,\cA=l(\gt g^*)$ was conjectured by Mishchenko and Fomenko \cite{mf}, and proved by Sadetov \cite{sad}. A clearer
treatment of this result is given by Bolsinov \cite{bol2}. Note that of course the image of a Poisson-commutative subalgebra $\cA\subset \cS(\gt g)$ remains Poisson-commutative
in $ {\mathcal S}(\gt g)/I$.
However, transcendence degree may sink far below $l(I)$. Our proof of Theorem~\ref{main2} follows the same strategy as the proofs of Sadetov and Bolsinov for ${\mathcal S}(\gt g)$. Note that in the general case our functions
$f_1, \ldots, f_{l(I)}\in  {\mathcal S}(\gt g)/I$ do not extend to Poisson-commuting functions in $\cS(\gt g)$.

First we prove Theorem~\ref{main2} in case of a reductive $\gt g$, see Section~\ref{red}. In the general case, we argue by induction on $\dim\gt g$, see Section~\ref{proof}.
We remark that the number $l(I)$ does not change under field extensions.   

Some applications of  Theorems~\ref{main-new} and \ref{main2}
to integrable Hamiltonian systems are discussed in
Section~\ref{appl}.

\section{Symplectic algebras and Hamiltonian actions}
\label{sh}

Consider a Poisson algebra ${\cP}$. Assume that ${\cP}$ has no zero-divisors  and that
${\rm tr.deg}\,{\cP}< \infty$.

For each subalgebra
 $C\subset{\cP}$,  let $C^{-1}{\cP}$ denote the localisation of
 ${\cP}$ by the subset of all non-zero elements of $C$. Clearly
$C^{-1}{\cP}$ is a subset of the field  ${\rm Quot}\,{\cP}$.
The Poisson structure uniquely  extends from ${\cP}$ to
${\rm Quot}\,{\cP}$ and for any multiplicative system
$S\subset{\cP}$ the localisation
${\cP}_S$ is a Poisson subalgebra of ${\rm Quot}\,{\cP}$. In particular, this is true for $C^{-1}{\cP}$. If $C\subset {\eus Z}{\cP}$, then
$C^{-1}{\cP}$ can be regarded as a Poisson algebra over the field
${\rm Quot}\, C$.

 \begin{df} A Poisson algebra ${\cP}$ is said to be {\it separable} if
${\rm tr.deg}\,{\eus Z}{\cP}+\rk{\cP}={\rm tr.deg}\,{\cP}$.
 \end{df}

Roughly speaking, ${\cP}$ is separable if generic symplectic leaves of
the underlying Poisson affine variety $X={\rm Spec}({\cP}\otimes_{\fk}\overline{\fk})$ are separated by the ``central" functions, elements of  ${\eus Z}{\cP}\otimes_{\fk}\overline{\fk}$.

If $\cP$ is a separable Poisson algebra, then
$({\eus Z}\cP)^{-1}\cP$ is a symplectic algebra over
${\rm Quot}\,{\eus Z}\cP$, see Definition~\ref{sympl}.

\begin{ex} Let $W$ be a finite-dimensional vector space over $\fk$ equipped with a non-degenerate skew-symmetric bilinear form $\omega$. Then $\omega$ defines a Poisson bracket on the symmetric algebra ${\mathcal S}(W)$ by the formula
$$
\{x,y\}:=\omega(x,y) \enskip \text{ for all $x,y\in W$.}
$$
This Poisson algebra $({\mathcal S}(W),\omega)$ is symplectic and
$V({\mathcal S}(W))={\rm Quot}\,{\mathcal S}(W)\otimes_{\fk} W$ with the same (extended) form
$\omega$.

The algebra ${\mathcal S}(W)$ has a natural grading, with grading components being
${\mathcal S}^k (W)$, and $\{{\mathcal S}^k(W),{\mathcal S}^l (W)\}\subset
 {\mathcal S}^{k+l-2}(W)$. Hence $\gt q:=\fk+W+{\mathcal S}^2(W)$ is a Lie subalgebra and $\gt n:=\fk+W$ is an ideal of $\gt q$. Note that $\gt n$ is a Heisenberg
algebra. The map
$$
\begin{array}{c}
\gt q\to {\rm Der}\,\gt n \\
\delta\mapsto \ad(\delta)_{|\gt n}\\
\end{array}
$$
is an epimorphism  of  Lie algebras 
with the
kernel $\fk$ and ${\mathcal S}^2(W)$ is mapped isomorphically onto the Lie algebra
$\gt{sp}(W)$ of the symplectic group $\Sp(W,\omega)$.
Note also that the centraliser of $\gt n$ in ${\mathcal S}(W)$ coincides with $\fk$.
\end{ex}

\begin{ex} Let $\gt q$ be a (finite-dimensional) Lie algebra over $\fk$. Then
${\mathcal S}(\gt q)$ is a Poisson algebra with  the usual Kirillov-Kostant  bracket. The corresponding symplectic vector space $V({\mathcal S}(\gt q))$ can be constructed as follows. Set $\mathbb K:={\rm Quot}\,\cS(\gt q)$. 
Recall that $\mathbb K$ is also a Poisson algebra. 
Set further $\widetilde{\omega}(\xi,\eta):=[\xi,\eta]$ for all $\xi,\eta\in\gt q$. Since $\widetilde{\omega}(\xi,\eta)\in \mathbb K$, this formula defines a skew-symmetric bilinear form on a $\mathbb K$-vector space $\widetilde{V}:=\gt q\otimes_{\fk}\mathbb K$.
(In a basis of $\gt q$, $\widetilde{\omega}$ is just the structural matrix.) Then
$V(\cS(\gt q))=\widetilde{V}/{\rm Ker}\,\widetilde{\omega}$. Let us say that $\gt q$ is {\it separable} if ${\rm tr.deg}\,{\eus Z}\cS(\gt q)={\rm tr.deg}\,{\eus Z}\mathbb K$. 
A Lie algebra $\gt q$ is separable if and only if $\cS(\gt q)$ is separable. In that case
$({\eus Z}\cS(\gt q))^{-1}\cS(\gt q)$ is a symplectic algebra over
${\rm Quot}\,{\eus Z}\cS(\gt q)$.
\end{ex}

The next two statements follow easily from the construction of $(V(\cP),\omega)$ and
Definition~\ref{sympl}.

\begin{prop}\label{prolog1} Let $\cP$ be a symplectic algebra and $\cA\subset\cP$ a Poisson subalgebra.
Let $U(\cA)\subset V(\cP)$ be a subspace spanned over ${\rm Quot}\,{\cP}$ by the derivations $\ad(\varphi)$ with
$\varphi\in\cA$.  Then

(1) $\dim_{{\rm Quot}\,{\cP}} U(\cA)={\rm tr.deg}\,\cA$;

(2)  $\rk \omega_{|U(\cA)}=\rk \cA$;

(3) $U(\cA)/{\rm Ker}\,w_{|U(\cA)}={\rm Quot}\,{\cP}\otimes_{{\rm Quot}\,{\mathcal A}}V(\cA)$.
\end{prop}

\begin{prop}\label{prolog2} Let $\cA, {\mathcal B}\subset\cP$ be two Poisson subalgebras
of a symplectic algebra $\cP$.
Then $\{\cA,{\mathcal B}\}=0$ if and only if the subspaces $U(\cA)$ and
$U({\mathcal B})$ are orthogonal with respect to $\omega$.
\end{prop}

Combining Propositions~\ref{prolog1} and \ref{prolog2}, we get that if
 $\{\cA,{\mathcal B}\}=0$, then
${\rm tr.deg}\,\cA+{\rm tr.deg}\,{\mathcal B}\le\rk\cP$.

From now on assume that $\cP$ is symplectic and that we have a Hamiltonian action of a Lie
algebra $\gt q$ on $\cP$, see Definition~\ref{ham}.
Set
$$
\cP^{\gt q}:=\{\varphi\in\cP \mid \{\gt q,\varphi\}=0\}={\eus Z}_{\cP}(\cP(\gt q)).
$$
As above, ${\rm tr.deg}\,\cP(\gt q)+{\rm tr.deg}\cP^{\gt q}\le\rk\cP$.
A Hamiltonian action is said to be {\it separable} if
${\rm tr.deg}\,\cP(\gt q)+{\rm tr.deg}\cP^{\gt q}=\rk\cP$.
%
%
It is possible to characterise separable Hamiltonian coisotropic actions.

\begin{prop}\label{prolog5} A separable Hamiltonian action
$\gt q\hookrightarrow\cP$
is coisotropic if and only if $\{\cP^{\gt q},\cP^{\gt q}\}=0$.
\end{prop}
\begin{proof} Recall that for a separable action, the orthogonal complement of
$U(\cP(\gt q))$ coincides with $U(\cP^{\gt q})$. Hence the action $\gt q\hookrightarrow\cP$ is coisotropic if and only if $U(\cP^{\gt q})\subset V(\cP)$ is an isotropic subspace. According to Proposition~\ref{prolog2}, this condition is equivalent to the Poisson-commutativity of $\cP^{\gt q}$.
\end{proof}

\begin{prop}\label{vyvod}
Theorem~\ref{main-new} follows from Theorem~\ref{main2}.
\end{prop}
\begin{proof}
The subalgebra $\cP(\gt q)\subset\cP$ is isomorphic to some Poisson quotient $\cS(\gt q)/I$, 
where $I\lhd\cS(\gt g)$ is a Poisson ideal. In particular, $I$ is $G$-invariant. 
Since $\cP$ is a domain, the algebra $\cS(\gt q)/I$ is also a domain.
By Theorem~\ref{main2}, $\cS(\gt q)/I$ contains a Poisson-commutative subalgebra $\cA$ with
${\rm tr.deg}\,\cA=l$, where
$$l= {\rm tr.deg}\,\cP(\gt q)- \frac{1}{2}\rk \cP(\gt q).$$
Combining Definition~\ref{cois} and Proposition~\ref{prolog1}, we see that
$$\rk \cP(\gt q)=\rk\cP-2(\rk\cP-{\rm tr.deg}\,\cP(\gt q))=2{\rm tr.deg}\,\cP(\gt q)-\rk\cP$$ and therefore
$l=\frac{1}{2}\rk\cP$.
\end{proof}

\section{Geometric realisation and Applications}\label{appl}  


Suppose that  $X$ is an irreducible affine  variety defined over 
$\fk$. 
Let $X(\overline\fk)$ denote  the set of its points over
the algebraic closure of $\fk$.  As usual 
$\fk[X]$ and $\fk(X):={\rm Quot}\,\fk[X]$ stand for the algebras
of regular and rational functions on $X$, respectively. 
Our convention is that $\overline\fk[X(\overline\fk)]=\fk[X]\otimes_{\fk}\overline\fk$. 
All subvarieties of $X$,  all
differential forms on $X$, and all morphisms of $X$ are
supposed to be defined over $\fk$.

Let $G$ be a connected linear algebraic group over $\fk$
with $\gt g=\Lie G$. 
An algebraic  action of $G$ on $X$ gives rise to a representation of $G$ (and of $\gt g$) on 
$\fk[X]$. 

\begin{df}[Geometric version of Definition~\ref{ham}]\label{Ham-geom} Suppose that $Y$ is an affine variety such that  $\fk[Y]$ is a Poisson algebra.
An algebraic action $G\times Y\to Y$ is said to be
{\it Hamiltonian} if there is a $G$-equivariant map, called
{\it the moment map}, $\mu\!: Y\to\gt g^*$ such that
$\mu^*:\, {\mathcal S}(\gt g)\to\fk[Y]$ is a Poisson homomorphism and 
$\{\mu^*(\xi),h\}=\xi{\cdot}h$ for all $\xi\in\gt g$, $h\in\fk[Y]$. 
\end{df}

Suppose that we have a Hamiltonian action $G\times Y\to Y$.
Then each function in $\mu^*({\mathcal S}(\gt g))$ is called
a {\it Noether integral} on $Y$. Their most important property
is given by the Noether theorem:
$\{\fk(Y)^G,\mu^*({\mathcal S}(\gt g))\}=0$.
The kernel of $\mu^*$ is a Poison ideal of $\cS(\gt g)$, say $I$, and therefore
$\cS(\gt g)/I$ is a Poisson quotient of $\cS(\gt g)$.

Let $I\lhd \cS(\gt g)$ be a $G$-invariant prime ideal. Being
$G$-invariant implies  that $\{\gt g,I\}\subset I$. In other words
$I$ is a Poisson ideal. Set $X= {\rm Spec}({\mathcal S}(\gt g)/I)$.
Then $\fk[X]$ is  a Poisson algebra and $X$ is a Poisson variety.
Set
$$
l(X):=(\dim X+{\rm tr.deg}\,\fk(X)^G)/2.
$$
It follows from Rosenlicht's theorem, that
$$l(X)=\frac{1}{2}(\dim X+(\dim X-\max_{\gamma\in X(\overline\fk)}\dim (\gt g\gamma))=
\dim X-
\frac{1}{2}\max_{\gamma\in X(\overline\fk)}\dim (\gt g\gamma),$$
where $X(\overline\fk)\subset(\gt g\otimes_{\fk}\overline\fk)^*$ and $\gt g\gamma=T_{\gamma} (G\gamma)$.
Let $R\subset\fk[X]$ be a Poisson commutative subalgebra.
Take $\gamma\in X(\overline\fk)$ such that the orbit
$G\gamma$ is of maximal possible dimension.
The subspace
$$\left<d_\gamma a\mid a\in R \right>\subset T^*_\gamma X(\overline\fk)\subset\gt g\otimes_{\fk}\overline\fk,$$
spanned over $\overline\fk$ by the differentials $d_\gamma a$, is isotropic with respect
to the symplectic form $\hat\gamma(x,y)=\gamma([x,y])$ (here $x,y\in\gt g\otimes_{\fk}\overline\fk$).
Hence the dimension of this subspace is less than or equal to $l(X)$ and
also ${\rm tr.deg}\,R\le l(X)$.
A family $\{f_1,\ldots, f_{l(X)}\}\subset\fk[X]$
is said to be {\it complete\/} if
$\{f_i,f_j\}=0$ for all $i,j$ and
$f_1,\ldots, f_{l(X)}$ are algebraically independent.

\vskip1ex

From now until the end of this section, assume that the geometric points of an irreducible affine variety 
$M$ form a dense subset of 
$M(\overline\fk)$. 
Suppose further  that there is a non-degenerate closed regular $2$-form
$\omega$ on the smooth locus of $M(\overline\fk)$. Then
$\omega$ induces a Poisson bracket on $\fk(M)$, and $\fk(M)$ is a symplectic Poisson algebra in the sense of
Definition~\ref{sympl}. The
variety $M$ is said to be {\it symplectic} if $\fk[M]$ is a
Poisson subalgebra of $\fk(M)$, i.e., if
$\{\fk[M],\fk[M]\}\subset\fk[M]$. In that case $\fk[M]$ is a symplectic algebra as well.  This is
always the case for normal affine varieties. Set
$2n:=\dim_\fk M={\rm tr.deg}\,\fk(M)$.
A family of functions $\{f_1,\ldots, f_n\}\subset\fk(M)$
such that $\{f_i,f_j\}=0$ for all $i$ and $j$
is said to be {\it complete} if the $f_i$'s are algebraically independent.
A complete family of functions generates a Poisson-commutative subalgebra
 $\cA\subset \fk(M)$ with ${\rm tr.deg}\,\cA=n$.

The simplest example of a symplectic variety is an
even-dimensional vector space $V$ equipped with a
non-degenerate skew-symmetric bilinear form $\omega$.
Each Lagrangian decomposition
$V=V_+\oplus V_-$ gives us a complete   family of
linear functions on $V$, namely, one has to take a basis of
$V_+^*$. Another familiar example is the cotangent bundle of a smooth irreducible
affine variety $Y$, $M=T^*Y$, equipped with
the canonical symplectic structure.
Here $\fk[Y]$ is a  Poisson commutative subalgebra of
$\fk[M]$. 
Since $\dim M =2\dim Y$, the subalgebra
$\fk[Y]$ contains a complete family of  functions on $M$.

It is a challenging open problem
to prove that for each affine symplectic variety $M$, the Poisson algebra $\fk[M]$
contains a complete family.

Suppose $h\in \fk[M]$. 
Let $\eta_h$ be the vector
field on the smooth locus of $M$ uniquely defined by the
formula $dh=\omega(\eta_h,\,.\,)$. Then $\eta_h$ defines
a Hamiltonian dynamical systems on $M$, and any function $f$
on $M$ such that $\{h,f\}=0$ is called a {\it first
integral} of this system. The intersection of the level
hypersurfaces of first integrals  is stable with respect to
the flow generated by $\eta_h$. Thus, to understand
dynamical properties of $\eta_h$, it is desirable to
construct as many independent first integrals as possible.
The triple  $(M, \omega, h)$ is said to be {\it completely
integrable} if there are algebraically independent first
integrals $f_1,\ldots,f_n$ such that $\{f_i,f_j\}=0$ and $2n=\dim M$.




Let $(M,\omega)$ be a symplectic affine variety and
$G\times M\to M$  a Hamiltonian algebraic action.
Write $M^{\rm reg}$ for the smooth locus of $M$.
For each $x\in M^{\rm reg}$, let
$(\gt gx)^{\perp}\subset T_xM$ denote the orthogonal complement of $\gt gx$ taken
with respect to $\omega$.
The function $$x\mapsto\dim(\gt gx\cap(\gt gx)^{\perp})$$ is
constant on a non-empty open subset $U\subset M^{\rm reg}$ and
its value $d$ on $U$ is called the  {\it defect} of the action $G\times
M\to M$ (see \cite[Chapter II, \S 3]{Vin}).

\begin{df}
The {\it corank} of $G\times M\to M$, denoted by
${\rm cork}\,M$, is defined by the formula
\[
   {\rm cork}\,M:=\min_{x\in M^{\rm reg}}\dim(\gt gx)^{\perp}-d \ .
\]
In other words, it equals  the rank of the form
$\omega|_{(\gt gx)^{\perp}}$ for generic $x\in M^{\rm reg}$.
A Hamiltonian action of $G$ on a symplectic variety $M$
is said to be {\it coisotropic} if ${\rm cork}\,M=0$, i.e., if
$(\gt gx)^{\perp}\subset \gt gx$ for generic $x\in M^{\rm reg}$.
\end{df}

\begin{thm}[A geometric version of Theorem~\ref{main-new}]      \label{main}
Let $G\times M\to M$ be a coisotropic  Hamiltonian action on a symplectic variety $M$ 
and let $\mu: M\to \gt g^*$ be the corresponding moment map.
Then there are functions $f_1,\ldots,f_n\in{\mathcal S}(\gt g)$,
where $n=\dim M/2$, such that
$\{\mu^*(f_1),\ldots,\mu^*(f_n)\}$ is a complete family
on $M$.
\end{thm}
\begin{proof}
Let $I\lhd \cS(\gt g)$ be the kernel of $\mu^*$. Then $I$ is a prime Poisson ideal of $\cS(\gt g)$. 
Set $X:= {\rm Spec}({\mathcal S}(\gt g)/I)$ and take 
$x\in M$. 
Using the fact that $M$ is a symplectic variety and the property 
$\{\mu(\xi),h\}=\xi{\cdot}h$ of the moment map, see Definition~\ref{Ham-geom},
one deduces that  the  kernel of $d\mu_x$  
 coincides with $(\gt gx)^{\perp}$, cf. \cite[Chapter~II, \S 2]{Vin}.
Therefore,
$\dim X=\max_{x\in M}\dim(\gt g x)$ and
$\max_{\gamma\in X(\overline\fk)}\dim_{\overline\fk}(\gt g\gamma)= \max_{x\in
M}\dim_{\fk}(\gt g x)-d$. Choose $x\in M$
such that $\dim(\gt gx)$ is maximal. Then
$$
l(X)=l(I)=\dim(\gt g x)-
\frac{1}{2}\left(\dim(\gt g x)-d\right)=
(\dim(\gt gx)+d)/2=
\left(\dim M-{\rm cork}\,M\right)/2.
$$
Clearly, $2l(X)=\dim M$ if and only if the action
$G\times M\to M$ is coisotropic.
By virtue of Theorem~\ref{main2},
there is a complete  family $\{f_i\}$ in ${\mathcal S}(\gt g)/I=\fk[X]$. Since ${\rm cork}\,M=0$ and $\mu^*$ is a Poisson homomorphism,
$\{\mu^*(f_i)\}$ is a complete family on $M$.
\end{proof}

\vskip0.3ex
\noindent
{\bf Corollary.} {\it
A Hamiltonian action $G\times M\to M$ is coisotropic
if and only if there is a complete family
of Noether integrals on $M$; or, equivalently, each
$G$-invariant Hamiltonian system on $M$ is
completely integrable in the class of Noether integrals.
}
\vskip0.4ex

\begin{thm}\label{int-2}
Let $G\times M\to M$ be a Hamiltonian action with
${\rm cork}\,M=2$. Then there is a complete
family in $\fk(M)$. If in addition
generic $G(\overline\fk)$-orbits on $M(\overline\fk)$ are separated by regular
invariants, then there is
a complete family in $\fk[M]$.
\end{thm}
\begin{proof}
Let $I\lhd \cS(\gt g)$ be the kernel of $\mu^*$. Set $X=: {\rm Spec}({\mathcal S}(\gt g)/I)$.  Then
$l(X)=\dim M/2-1$. By Theorem~\ref{main2},
there are functions $f_1,\ldots,f_{l(X)}\in{\mathcal S}(\gt g)$
such that their restrictions to $X$ form a complete
family. Set $R:=\mu^*({\mathcal S}(\gt g))$.
Let $\left<d_x R\right>$ be the subspace of $T^*_x M$
spanned over $\fk$ by all differentials $d_x f$ with $f\in R$.
Since $(\gt gx)^{\perp}$ is the kernel of $d\mu_x$,
we have $\left<d_x R\right>={\rm Ann}\,((\gt gx)^{\perp})$.
By Rosenlicht's theorem,
generic $G(\overline\fk)$-orbits on $M(\overline\fk)$ are separated by rational invariants.
Therefore,  $\left<d_x(\fk(M)^G)\right>={\rm Ann}\,(\gt gx)$ for generic $x\in M$. Since the action
$G\times M\to M$ is not coisotropic,
$(\gt gx)^{\perp}\not\subset \gt gx$ and there is at least
one $h\in\fk(M)^G$ such that functions $\{h,\mu^*(f_1),\ldots,\mu^*(f_{l(X)})\}$
are algebraically independent.
Recall that $\{\fk(M)^G,R\}=0$. Thus,
$\{h,\mu^*(f_1),\ldots,\mu^*(f_{l(X)})\}$ is a complete family
on $M$. If
generic $G(\overline\fk)$-orbits on $M(\overline\fk)$ are separated by regular
invariants, then $\fk(M)^G={\rm Quot}\,\fk[M]^G$ and
we can choose $h$ in $\fk[M]^G$.
\end{proof}

Let us say a few words about
cotangent bundles. It was already mentioned that
a complete family always exists here.
But the construction
of Theorems~\ref{main2} and \ref{main} provides other
examples of complete families, which can be useful
for other Hamiltonian systems.

Suppose that $M=T^*X$,
where $X$ is a $G$-variety. Then $M$ possesses a canonical
$G$-invariant symplectic structure such that the action
of $G$ is Hamiltonian.  If
the action $G\times M\to M$ is coisotropic,
then $X$ has an open $G$-orbit \cite{gu}.
For reductive $G$ one can say more.

Suppose $\fk$ is algebraically closed, $G$ is reductive, and $X$ is smooth.
By a result of Knop~\cite[Sections~6\&7]{knop},
the action of $G$ on $T^*X$ is a coisotropic if and only if
a Borel subgroup $B$ of $G$ has on open orbit on $X$.
Normal varieties having an open $B$-orbit are
said to be {\it spherical}.
It was known before that if $X$ is spherical and $X=G/H$, where
$H$ is a reductive subgroup of $G$, then
each $G$-invariant Hamiltonian
system on $T^*X$ is integrable within
the class of Noether integrals, see \cite{gu,mi,alan}.
Here we 
lift the assumption that $H$ is reductive.
Smooth affine spherical varieties are classified
(under mild technical constraints) in \cite{b-kn}.
It would be interesting to study complete families
on their cotangent bundles.

By the same result of Knop \cite{knop}, the action of $G$ on $T^*X$ is of corank $2$ if and only if
${\rm tr.deg}\,\fk(X)^B=1$, i.e., $X$ has {\it complexity} $1$.
 Theorem~\ref{int-2} provides also (hopefully) interesting completely integrable systems for these
 cotangent bundles.

Other well-studied coisotropic actions on cotangent bundles
are related to  {\it Gelfand pairs}.
Suppose that $\fk=\rl$ and
$M=T^*X$, where $X=G/K$ is a Riemannian homogeneous space.
Then  $X$ is called
{\it commutative} or the pair
$(G, K)$ is called a {\it Gelfand pair} if the action
$G\times M\to M$ is coisotropic. Gelfand pairs can be
characterised by the following equivalent conditions.
\begin{itemize}
\item[({\sf i})] \ The algebra ${\mathcal D}(X)^G$ of $G$-invariant differential
operators on $X$
is commutative.
\item[({\sf ii})] \ The algebra of $K$-invariant measures on $X$ with
compact support is commutative with respect to convolution.
\item[({\sf iii})] \ The representation of $G$ on
$L^2(X)$ has a simple spectrum.
\end{itemize}
Theorem~\ref{main} and its corollary provide
 two more equivalent conditions.
\begin{itemize}
\item[({\sf iv})] \ There is a complete family
of Noether integrals on $T^*X$.
\item[({\sf v})] \ Each $G$-invariant Hamiltonian system on $M$ is
completely integrable in the class of Noether integrals.
\end{itemize}

According to \cite{Vin}, if $G/K$ is a Gelfand pair and
$G=L\ltimes N$ is a Levi decomposition of $G$ such that
$K\subset L$, then $\rl[\gt n]^L=\rl[\gt n]^K$ and
$\gt n$ is at most two-step nilpotent. These conditions
guarantee that the construction of a complete family on
$\mu(M)$ would have at most three induction steps.
Thus, one can hope for explicit formulas for
our commuting families and applications to
physical problems. Gelfand pairs are partly classified
in \cite{Vin2,tg} and completely in  \cite{dis}.

\section{The reductive case}   \label{red}

In this section, $G$ is a connected reductive algebraic group.
Here one can apply a very powerful tool, the
so called ``argument shift method''. It was used by
Manakov \cite{man}, Mishchenko and Fomenko \cite{mf3}, and
Bolsinov \cite{bol}
in constructions of complete families on $\gt g^*$ and
coadjoint $G$-orbits. The reader is referred to \cite[Chapter 4]{tr-f1}
for a thorough exposition and
historical remarks. Let us briefly outline this method.

Let $r$ be the rank of $\gt g$.
Choose any set $F_1,\ldots,F_r$
of free generators of $\fk[\gt g^*]^G$.
For any $a\in\gt g^*$, let  ${\mathcal F}_a$ denote the finite set
$$
\{F_i, \partial_a F_i, \partial_a^2 F_i, \ldots ,\partial_a^{k(i)} F_i\mid i=1,\ldots,r, \ k(i)=\deg F_i-1\}\subset{\mathcal S}(\gt g).$$
Then $\{{\mathcal F}_a,{\mathcal F}_a\}=0$,
see e.g. \cite[Sections~1.12,\,1.13]{per}.
Here we should mention that
this fact is stated in \cite{per} for $\fk=\cp$,
but the proofs are valid over all fields of characteristic zero.

Recall that the {\it index of a Lie algebra} $\gt q$
is the minimum of dimensions of
stabilisers $\gt q_\xi$ over all covectors $\xi\in\gt q^*$,
i.e., $\ind\gt q=\min\limits_{\xi\in\gt q^*}\dim\gt q_\xi$. Note that
$\ind\gt q={\rm tr.deg}\,\fk(\gt q^*)^{\gt q}$ and that $\dim\gt q-\ind\gt q$ is the rank of the Poisson
algebra $\cS(\gt q)$ as defined in the Introduction.

\begin{prop}{\cite[Theorem 2]{bol}}   \label{bol}
Suppose that  $\gt g$ is a {\sl complex}
reductive Lie algebra and $\xi\in\gt g^*$.
Then there is $a\in\gt g^*$ such that the
restriction of ${\mathcal F}_a$ to the
coadjoint orbit $G\xi$
contains $\frac{1}{2}\dim(G\xi)$
algebraically independent functions
if and only if $\ind\gt g_\xi=\ind\gt g$.
\end{prop}

The proof of Theorem~2 in \cite{bol} uses only linear algebra and can be repeated for 
any algebraically closed field of characteristic zero. We are going to use the result also for $\fk\ne \mathbb C$.


\begin{prop}\label{Elashvili}
If $\gt g$ is reductive and $\xi\in\gt g^*$, then
$\ind\gt g_\xi=\ind\gt g$.
\end{prop}

The statement of Proposition~\ref{Elashvili} is known as Elashvili's conjecture.
For the classical Lie algebras,  it is proved in \cite{ya1}
under the assumption that $\chr\fk$ is good for $\gt g$.
W.~de Graaf used a computer
program to verify the conjecture for the exceptional Lie algebras,
see \cite{graaf}.
An almost conceptual  proof of Elashvili's conjecture is
given in \cite{Ch-M}. (The authors still have to rely on  computer calculations for a few orbits.)

Let $\hat V_{a,\xi}\subset T^*_\xi(\gt g^*)$  be the $\fk$-linear span of 
 the differentials $\{d_\xi F\mid F\in{\mathcal F}_a\}$
and let $V_{a,\xi}$ be the restriction of
$\hat V_{a,\xi}$ to $T_\xi(G\xi)=\gt g\xi$.
Since the orbit $G\xi$ is a symplectic variety and
the subspace $V_{a,\xi}$ is isotropic, we get
$2\dim V_{a,\xi}\le\dim(G\xi)$.
The restriction of ${\mathcal F}_a$ to
$G\xi$ contains a complete family if and only if
there is $a'\in Ga$ such that
$2\dim V_{a',\xi}=\dim G\xi$.

Combining Propositions~\ref{bol} and \ref{Elashvili},
we obtain the following assertion.

\begin{prop}    \label{point}
Suppose that $\fk=\overline{\fk}$.
Then for each $\xi\in\gt g^*$, there is $a\in\gt g^*$ such that
$2\dim V_{a,\xi}=\dim G\xi$.
\end{prop}

%

\begin{proof}[Proof of Theorem~\ref{main2} in the reductive case.]
Let $I$ be a prime Poisson ideal of ${\mathcal S}(\gt g)$ and
$X(\overline\fk)$  a closed subvariety of
$(\gt g\otimes_{\fk}\overline\fk)^*$ defined by $I$.

Choose a set of homogeneous generators
$\{F_1,\ldots,F_r\}\subset\fk[\gt g^*]^G$.
Let $\hat{F_i}$ denote the restriction of $F_i$ to $X$.
Each fibre of the quotient morphism $X(\overline\fk)\to X(\overline\fk)/\!\!/G(\overline\fk)$ contains finitely many
$G$-orbits.
Hence for generic $\xi\in X(\overline\fk)$
the differentials $\{d_\xi \hat{F_i}\mid i=1,\ldots,r\}$
generate a subspace of dimension $m:=\dim X-\dim(G\xi)$.
According to Proposition~\ref{point}, there is
an element $a\in(\gt g\otimes_{\fk}\overline\fk)^*$ such that
the restriction of ${\mathcal F}_a$ to $G(\overline\fk)\xi$
contains a complete family, i.e., $2\dim V_{a,\xi}=
\dim(G(\overline\fk)\xi)$. There is an open subset of such elements.
In particular, we may (and will) assume that $a\in\gt g^*$. Then
${\mathcal F}_a$ is a subset of ${\mathcal S}(\gt g)$.
Each differential $d_\xi\hat{F_i}$ is zero on $\gt g\xi$. Therefore
$$
\dim\left<d_\xi f\mid f\in{{\mathcal F}_a}/({\mathcal F}_a\cap I) \right>=
 m+\dim(G\xi)/2=(m+\dim X)/2=l(X)
$$
and the restriction of ${\mathcal F}_a$ to $X$ contains a complete family.
\end{proof}


\section{Auxiliary results}\label{dop}

In this section, we collect several 
facts concerning structural properties of algebraic Lie
algebras. They will be used in the proof of the main theorem.

Recall that  a $(2n{+}1)$-dimensional Heisenberg Lie algebra over $\fk$
is a Lie algebra $\gt h$ with  a basis
$\{x_1,\ldots,x_n,y_1,\ldots,y_n,z\}$  such that $n\ge 1$,
$[x_i,x_j]=[y_i,y_j]=0$, $[\gt h,z]=0$, and $[x_i,y_j]=\delta_{ij}z$.
Recall also that a Lie ideal $\gt a\lhd\gt q$ is said to be a {\it characteristic ideal} if
it is stable under all automorphisms of the Lie algebra $\gt q$.

\begin{lm}\label{nil}
Suppose that $\gt n$ is a nilpotent Lie algebra
such that each commutative
characteristic ideal of $\gt n$ is one-dimensional.
Then $\gt n$ is a Heisenberg algebra.
\end{lm}

\begin{proof}
Let $\gt z$ be the centre of $\gt n$.
Then $\dim\gt z=1$. Consider  
the upper central series of $\gt n$
$$
\gt z=\gt n_0\subset\gt n_1\subset\gt n_2\subset\dots\subset\gt n_{k-1}
\subset\gt n_k=\gt n,
$$
i.e., $\gt n_i/\gt n_{i-1}$ is the centre
of $\gt n/\gt n_{i-1}$. The centre of $\gt n_1$
is a commutative characteristic ideal of $\gt n$.
Hence, it is one-dimensional and coincides with
$\gt z$. Therefore  $\gt n_1$ is a Heisenberg algebra.
Let $\gt z_{\gt n}(\gt n_1)$ be the centraliser
of $\gt n_1$ in $\gt n$. Clearly, $\gt z_{\gt n}(\gt n_1)$
is an ideal in $\gt n$ and $\gt n_1\cap\gt z_{\gt n}(\gt n_1)=\gt z$.
We claim that $\gt n=\gt n_1+\gt z_{\gt n}(\gt n_1)$.
Indeed, let $\xi\in\gt n$. Then
$[\xi,\gt n_1]\subset\gt n_0$ and
there is an element $\xi_0\in\gt n_1$ such that
$[\xi-\xi_0,\gt n_1]=0$.

Let $\gt z_0$ be the centre
of $\gt z_{\gt n}(\gt n_1)/\gt z$.
Since  $\gt n/\gt z=(\gt n_1/\gt z)\oplus
(\gt z_{\gt n}(\gt n_1)/\gt z)$ is the direct sum of
two ideals,
$\gt z_0$ lies in the centre of
$\gt n/\gt z$. Thus,
$\gt z_0\subset(\gt n_1/\gt z)$ and
$\gt z_0=0$. Since
$\gt z_{\gt n}(\gt n_1)/\gt z$ is a nilpotent
Lie algebra, we have $\gt z_{\gt n}(\gt n_1)/\gt z=0$,
and $\gt n=\gt n_1$ is a Heisenberg algebra.
\end{proof}

Let $N$ be the unipotent radical of an affine algebraic group $G$.
Set $\gt n:=\Lie N$. For any action
$P\times Y\to Y$ let $Y/P$ stand for the
set of $P$-orbits on $Y$.

\begin{lm}\label{heis} Suppose that $\gt n$
is a Heisenberg Lie algebra  and the centre 
$\gt z$ of $\gt n$ lies in the centre  of $\gt g$.
Given a non-zero $\alpha\in\gt z^*$,
set $Y_\alpha:=\{\gamma\in\gt g^*|\, \gamma_{|\gt
z}=\alpha\}$.
Then $Y_\alpha/N={\rm Spec}\,\fk[Y_\alpha]^N$;
the natural action of $G/N$ on $Y_\alpha/N$ is Hamiltonian in the sense of Definition~\ref{Ham-geom}
and the moment map
$\mu:\,Y_\alpha/N\to(\gt g/\gt n)^*$ is
a $G$-isomorphism.
\end{lm}

\begin{proof}
Choose a Levi decomposition $G=L\ltimes N$ and
let $V$ be an $L$-invariant complement of $\gt z$ in $\gt n$.
Set $S_\alpha:=\{\gamma\in\gt g^*|\, \gamma(V)=0, \gamma_{|\gt
z}=\alpha\}$. In other words,
$S_\alpha=(\gt g/\gt n)^*+\widetilde{\alpha}$,
where $\widetilde{\alpha}\in\gt g^*$,
$\widetilde{\alpha}(V)=0$,
$\widetilde{\alpha}(\gt l)=0$, and
$\widetilde{\alpha}_{|\gt z}=\alpha$.
Clearly $S_\alpha\subset Y_\alpha$.
Each point $\gamma\in Y_\alpha$ can be uniquely
presented as a sum
$$
\begin{array}{c}
\gamma=\beta+\ad^*(\eta)\um
\widetilde\alpha+\widetilde\alpha,
\mbox{ where } \beta(\gt n)=0 \mbox{ and } \eta\in V. \\
\mbox{Thus }\
N\gamma\cap S_\alpha=\gamma-\ad^*(\eta)\cdot\gamma
+\frac{1}{2}(\ad^*(\eta))^2\um \gamma=\{pt\}.\\
\end{array}
$$
We obtain the isomorphism $\mu:\,\, Y_\alpha/N\to (\gt g/\gt n)^*$, where
$\mu(N\gamma)$ is the unique point in $N\gamma\cap S_\alpha$.
Therefore $Y_\alpha/N$ is an algebraic variety (an affine
space) and $\fk[Y_\alpha/N]=\fk[Y_\alpha]^N$. For the rest
of the proof,   we fix the isomorphism $(\gt g/\gt
n)^*\cong\gt l^*$ given by the Levi decomposition $\gt g=\gt
l\oplus\gt n$ and the induced isomorphisms
$S_\alpha\cong(\gt l^*+\widetilde\alpha)\cong\gt l^*$, where the last one is given by choosing $\widetilde\alpha$ as the origin.

For each $\gamma\in S_\alpha$ and $l\in L$, we have
$\mu(l\cdot N\gamma)=l\gamma-\widetilde{\alpha}=l(\gamma-\widetilde{\alpha})=l\cdot\mu(N\gamma)$.
This shows that $\mu$ is $G$-equivariant.

It  remains to  prove that $\mu^*$ is a homomorphism of the Poisson
algebras ${\mathcal S}(\gt g/\gt n)$ and $\mathbb F[Y_\alpha]^N$, i.e.,
to show that
$\{\mu^*(f_1),\mu^*(f_2)\}=\mu^*(\{f_1,f_2\})$ for all $f_1,f_2\in S(\gt g/\gt n)$.

Let $\gamma\in S_\alpha$. 
The identification ${\mathcal S}(\gt g/\gt n)\cong {\mathcal S}(\gt l)\subset{\mathcal S}(\gt g)$ gives
us that 
$$
\{f_1,f_2\}(\gamma)=\{f_1,f_2\}(\gamma-
\widetilde\alpha)=\{f_1,f_2\}(\mu(N\gamma))=
\mu^*(\{f_1,f_2\})(N\gamma).
$$
The last step is to prove that  $\{\mu^*(f_1),\mu^*(f_2)\}(N\gamma)=\{f_1,f_2\}(\gamma)$.
It is well-known that
$G\gamma$ is a symplectic leaf
of $Y_\alpha$ and $\gt g^*$.
Also $L\gamma$ is a symplectic leaf of $S_\alpha$.
We have
$$T_\gamma(G\gamma)=T_\gamma(L\gamma)\oplus
T_\gamma(N\gamma), \
\text{ where T$_\gamma(L\gamma)=\gt l\gamma$ \ and
\
$T_\gamma(N\gamma)=\gt n\gamma$ \ are orthogonal,} $$
additionally
 $T_\gamma(L\gamma)\subset T_\gamma S_\alpha$.
Let $F_i$ be the restriction of $\mu^*(f_i)$, which is regarded now as an $N$-invariant function on $Y_\alpha$,
to $G\gamma$. Since $G\gamma$ is a symplectic leaf of $\gt g^*$, we have
$$
\{\mu^*(f_1),\mu^*(f_2)\}(N\gamma)=\{F_1,F_2\}(\gamma).
$$
Clearly, the functions $F_1$ and $F_2$ are $N$-invariant, hence
$d_{\gamma}F_i({\gt n\gamma})=0$.
Thus
$\{F_1,F_2\}(\gamma)=\{F_1|_{L\gamma},F_2|_{L\gamma}\}(\gamma)=
\{f_1,f_2\}(\gamma)$ and we are done.
\end{proof}

\cled\
{\it In the setting of Lemma~\ref{heis}, we have
$$
(\fk[\gt g^*][1/z])^N\cong {\mathcal S}(\gt g/\gt n)\otimes_{\fk}\fk[z,1/z]\subset
{\mathcal S}(\gt g/\gt n)\otimes_{\fk}\fk(\gt z^*),
$$
where $z$ is a non-zero element of $\gt z$.
Moreover,
if $X(\overline\fk)\subset(\gt g\otimes_{\fk}\overline\fk)^*$ is a closed $G(\overline\fk)$-invariant subset
defined over $\fk$ and such that  $z_{|X(\overline\fk)}\ne 0$, then
$(\fk[X][1/z])^N$ is a Poisson quotient of
${\mathcal S}(\gt g/\gt n)\otimes_{\fk}\fk[z,1/z]$.
}

\begin{proof}
Suppose first that $\fk=\overline\fk$. Then $X(\overline\fk)=X$.
Set $X_\alpha:=X\cap Y_\alpha$.
Then $S_\alpha\cap X_\alpha$ defines a section of
$X_\alpha/N$, i.e.,
$X_\alpha/N\cong X_\alpha\cap S_\alpha\cong
{\gt x}_\alpha$, where ${\gt x}_\alpha\subset(\gt g/\gt n)^*$
is a $G$-invariant (Poisson) subvariety.
Therefore $(\fk[X][1/z])^N$ is a Poisson quotient of
${\mathcal S}(\gt g/\gt n)\otimes_{\fk}\fk[z,1/z]$.

Consider now the general case.
The Galois group ${\rm Gal}_\fk(\overline\fk)$ of the field
extension $\fk\subset\overline\fk$ acts on $(\overline\fk[X][1/z])^N$ and on 
${\mathcal S}(\gt g/\gt n)\otimes_{\fk}\overline\fk[z,1/z]$.
Taking its fixed points on both sides, we see that the statement holds.
\end{proof}

\begin{rmk} From Lemma~\ref{heis} 
one can deduce that
$\fk(\gt g^*)^G=\fk((\gt g/\gt n)^*)^{G/N}\otimes_{\fk}\fk(\gt z^*)$.
In particular, in this case $\fk(\gt g^*)^G$ is a
rational field.
\end{rmk}

Let $H\lhd N$ be a connected
commutative normal subgroup of $G$ with $\Lie H=\gt h$.

\begin{lm}\label{com}
Fix $\alpha\in\gt h^*$ and let
$Y_\alpha$
be the preimage of $\alpha$ under the natural
restriction $\gt g^*\to\gt h^*$. Then
$Y_\alpha/H={\rm Spec}\,\fk[Y_\alpha]^H$
and the restriction map
$\pi_\alpha:\,\,Y_\alpha\to(\gt g_\alpha)^*$ defines an isomorphism
$Y_\alpha/H\cong(\gt g_\alpha/\gt h)^*\times\{\alpha\}$.
\end{lm}
\begin{proof} Let $\gamma\in\gt g_\alpha$, $\xi\in\gt h$, $\eta\in\gt g$.
Then
$$
(\xi\um \gamma)(\eta)=\gamma([\eta,
\xi])=\alpha([\eta, \xi])=-(\eta\um\alpha)(\xi).
$$
Note that $\xi{\cdot}(\xi{\cdot}\gamma)=0$ and
therefore $H\gamma=\gamma+\gt h\gamma=\gamma+(\gt
g/\gt g_\alpha)^*$. Each non-zero fibre of the natural
$G_\alpha$-equivariant restriction $\pi_\alpha:\,
Y_\alpha\to(\gt g_\alpha)^*$ is exactly one $H$-orbit. Let
us fix a decomposition $\gt g=\gt g_\alpha\oplus\gt m$.
Choose any
$\widetilde\alpha\in\gt g^*$ such that
 $\widetilde\alpha(\gt m)=0$ and $\widetilde\alpha|_{\gt h}=\alpha$. Then
$Y_\alpha=(\gt g/\gt h)^*+\widetilde{\alpha}$ and
$\pi_\alpha(Y_\alpha)\cong
(\gt g_\alpha/\gt h)^*\times\{\widetilde\alpha\}\cong
 (\gt g_\alpha/\gt h)^*\times\{\alpha\}$.
\end{proof}

Until Lemma~\ref{non-closed}, we assume that $\fk=\overline\fk$.
Suppose that  $X\subset\gt g^*$ is a closed $G$-invariant subset.
Let $\gt x_{\gt h}\subset\gt h^*$ denote the image of
$X$ under the restriction $\gt g^*\to \gt h^*$.
Set $\mathbb K:=\fk(\gt x_{\gt h})$ and
\begin{gather}
\hat{\gt g}:=\{\xi\in \gt g\otimes_{\fk}\mathbb K \mid
    [\xi,\gt h](\gt x_{\gt h})=0\}, \label{g1}\\
\hat{\gt h}:=\{\xi\in\gt h\otimes_{\fk}\mathbb K\mid
\alpha(\xi(\alpha))=0 \mbox{ for each } \alpha\in\gt x_{\gt h}
   \ \text{ such that $\xi(\alpha)$ is defined}\}. \label{g2}
\end{gather}
Then $\hat{\gt g}$ is the Lie algebra
of all rational maps $\xi:\,{\gt x_{\gt h}}\to\gt g$ such that
$\xi(\alpha)\in\gt g_\alpha$ whenever $\xi(\alpha)$ is defined.

Since $\gt h$ is a commutative ideal of $\gt g$,  we have $\gt h\otimes_{\fk}\mathbb K\lhd\hat{\gt g}$.
Moreover,  $\hat{\gt h}$ is also an ideal of $\hat{\gt g}$.
The main object of our interest is the quotient Lie algebra
$\widetilde{\gt g}:=\hat{\gt g}/\hat{\gt h}$.
Another way to define this Lie algebra is to say that
$\widetilde{\gt g}:=\{\xi\in \gt g\otimes_{\gt h}\mathbb K
\mid [\xi,\gt h](\gt x_{\gt h})=0\}$.

Set 
${\mathcal A}:=(\fk[X]\otimes_{\fk[\gt x_{\gt h}]}
\mathbb K)^H=\fk[X]^H\otimes_{\fk[\gt x_{\gt h}]}\mathbb K$.
Then the algebra ${\mathcal A}$
carries a natural Poisson structure induced from $\fk[X]$.

\begin{lm}\label{com2}
Suppose that $\fk=\overline\fk$. Then
${\mathcal A}$ is a Poisson quotient of
${\mathcal S}(\widetilde{\gt g})$.
\end{lm}
\begin{proof}
The elements of
${\mathcal A}$ and
${\mathcal S}(\widetilde{\gt g})$ are linear combinations of  rational functions on
$\gt x_{\gt h}$ with coefficients from $\mathbb F[X]^H$ or
${\mathcal S}(\gt g)$, respectively. Thus, it suffices to verify the
claim at generic $\alpha\in\gt x_{{\gt h}}$.

Fix a vector space decomposition
$\gt g=\gt g_\alpha\oplus\gt m$ and let
$s:\,\{\alpha\}\times(\gt g_\alpha/\gt h)^*\to {\rm Ann}(\gt m)\subset Y_\alpha$
be the corresponding section of $\pi_\alpha$.
Then $S_\alpha:={\rm Im}\,s$ is a closed subset of $Y_\alpha$
and
by Lemma~\ref{com},
$S_\alpha\cap X\cong\pi_\alpha(Y_\alpha\cap X)\cong (X\cap Y_\alpha)/H$.
Let ${\mathcal A}_{\alpha}\subset\cA$ be the subset of elements  that are defined
at $\alpha$. Then  for generic $\alpha\in\gt x_\gt h$, we have a surjective map
$$
\epsilon_\alpha:\enskip {\mathcal A}_{\alpha}\to
  \fk[Y_\alpha\cap X]^H\cong\fk[S_\alpha\cap X].
$$
%
At the same time,
$\hat{\gt g}(\alpha):=\{\xi(\alpha)\mid\xi\in\hat{\gt g}, \xi(\alpha) \text{ is defined}\}=
\gt g_\alpha$ for generic $\alpha\in\gt x_{\gt h}$.
The algebra $\widetilde{\gt g}(\alpha):=\{\xi(\alpha)\mid\xi\in\widetilde{\gt g}, \xi(\alpha) \text{ is defined}\}$ 
is a $1$-dimensional central extension of $\gt g_\alpha/\gt h$. We have
$(\gt g_\alpha/\gt h)\oplus\fk w$ with the Lie bracket
$$
[\xi+\gt h,\eta+\gt h]:=([\xi,\eta]+\gt h)+\widetilde{\alpha}[\xi,\eta]w
\enskip
\text{ for all } \xi,\eta\in\gt g_\alpha,
$$
where $\widetilde{\alpha}\in\gt g_\alpha^*$ is
a linear function such that $\widetilde{\alpha}_{|\gt h}=\alpha$.
Hence
$(\widetilde{\gt g}(\alpha))^*=
  (\gt g_\alpha/\gt h)^*\times\fk\widetilde\alpha\,$ and
$S_\alpha\cap X$ is a closed subset of $(\widetilde{\gt g}(\alpha))^*$.

Therefore
${\mathcal A}$ is a quotient of
${\mathcal S}(\widetilde{\gt g})$. Since
the Poisson structure on ${\mathcal A}$ is induced
from $\fk[X]$ and $X$ is a Poisson subvariety of $\gt g^*$,
it is indeed a Poisson quotient.
\end{proof}

\begin{rmk} Informally speaking,
${\mathcal A}$ is the algebra of functions on the
set $\widetilde X$ of all
rational morphisms
$\psi:\,\gt x_{\gt h}\to X$ such that
$\psi(\alpha)\in (X\cap S_\alpha)$.
Here  $\widetilde X$ is also a set of the $H$-invariant
rational morphisms $\psi':\,\gt x_{\gt h}\to X$ such that
$\psi'(\alpha)\in (X\cap Y_\alpha)$.
\end{rmk}

If $\fk\ne\overline{\fk}$, then it is better to work with ideals. 
Let $I\lhd{\mathcal S}(\gt g)$ be 
a $G$-invariant prime ideal. Set $I_0=I\cap{\mathcal S}(\gt h)$ and let
$\gt x_{\gt h}\subset \gt h^*$ be the subvariety defined by $I_0$. 
Now $\mathbb K={\rm Quot}\,{\mathcal S}(\gt h)/I_0$ and
$$\widetilde{\gt g}:=\{\xi\in\gt g\otimes_{\gt h}\mathbb K\mid \{\xi,\gt h\}\subset I_0\otimes\mathbb K\}.$$
Finally set $\widetilde{\cP}:=({\mathcal S}(\gt g)/I)^{\gt h}\otimes_{\fk[\gt x_{\gt h}]}\mathbb K$.

\begin{lm}\label{non-closed}
Let $\fk$ be any field of characteristic zero. Then 
$\widetilde{\cP}$ is a Poisson
quotient of ${\mathcal S}(\widetilde{\gt g})$.
\end{lm}
\begin{proof}
In case  $\fk=\overline{\fk}$, $\widetilde{\gt g}$ coincides
with the quotient $\hat{\gt g}/\hat{\gt h}$, where
$\hat{\gt g}$ and $\hat{\gt h}$ are defined by Formulas~(\ref{g1})  and (\ref{g2}).
In the general case, we have
$\widetilde{\gt g}\otimes_{\fk}\overline\fk=
\widehat{\gt g\otimes_{\fk}\overline\fk}/\widehat{\gt h\otimes_{\fk}\overline\fk}$.
By Lemma~\ref{com2}, 
$\widetilde{\cP}\otimes_{\fk}\overline{\fk}$ is a Poisson quotient of 
${\mathcal S}(\widetilde{\gt g})\otimes_{\fk}\overline{\fk}$. 
The Galois group ${\rm Gal}_\fk(\overline\fk)$ of the field
extension $\fk\subset\overline\fk$ acts  on both these Poisson algebras.
By taking fixed points of  ${\rm Gal}_\fk(\overline\fk)$,
we conclude that
$\widetilde{\cP}$ is a Poisson
quotient of ${\mathcal S}(\widetilde{\gt g})$.
\end{proof}


\section{Inductive argument}\label{proof}

Let $I\lhd{\mathcal S}(\gt g)$ be a $G$-invariant (i.e., Poisson)
prime ideal, set $X={\rm Spec}({\mathcal S}(\gt g)/I)$. 
Then ${\cP}:=\fk[X]={\mathcal S}(\gt g)/I$ is a Poisson algebra. In this section,
we construct a complete  family in $\fk[X]$.

\begin{proof}[Proof of Theorem~\ref{main2}]
Set $n:=\dim X={\rm tr.deg\,}{\cP}$, $m:=\dim X-\rk{\cP}$.
Then $n-m$ is the dimension of a generic $G(\overline\fk)$-orbit on $X(\overline\fk)$, and $l=l(X)=(n+m)/2$.
The task is to construct
$l$ functions $f_i\in{\mathcal S}(\gt g)$ such that
$\{f_i,f_j\}\in I$ and their restrictions to $X(\overline\fk)$ are
algebraically independent. We argue by induction on $\dim\gt g$.
At first it is assumed that $\gt g$ is algebraic. 
The case of a non-algebraic Lie algebra $\gt g$ is treated at the very end.


 $\bullet$ \ In case of a reductive $G$, Theorem~\ref{main2} is proved in Section~\ref{red}.
Assume therefore that $G$ is not reductive. If $I$
contains a non-trivial ideal $\gt c\lhd\gt g$, then $\fk[X]$ is Poisson quotient of
$\cS(\gt g/\gt c)$  and we can replace $\gt g$ by $\gt g/\gt c$ without loss of
generality.
Below we assume that
$I\cap\gt g=0$. 
Let $\gt n$ be the nilpotent radical of $\gt g$ and
$N\subset G$ the connected subgroup
with $\Lie N= \gt n$.

$\bullet$ \ Suppose that $\gt n$ is a Heisenberg Lie algebra and
$\gt z=[\gt n,\gt n]$ is
a central subalgebra of $\gt g$. Then Lemma~\ref{heis} applies.  Let $z\in\gt z$ be a non-zero element.
Since $\gt z\not\subset I$, we have
$z_{|X(\overline\fk)}\ne 0$. Set
$\widetilde{\cP}:=({\cP}[1/z])^N$.
By Lemma~\ref{heis},
$\widetilde{\cP}$ is a Poisson
quotient of ${\mathcal S}(\gt g/\gt n)\otimes_{\fk}\fk[z,1/z]$.
The Lie
algebra $\gt g/\gt n$ is reductive, therefore there
are pairwise commuting functions
$f_1,\ldots, f_k\in{\mathcal S}(\gt g/\gt n)\otimes\fk(\gt z^*)$
such that
their images form a complete family in $\widetilde{\cP}$.
After multiplying by a
common denominator, we may assume that
each $f_i$ lies in ${\mathcal S}(\gt g/\gt n)$.

Choose a decomposition $\gt n=V_+\oplus V_-\oplus\gt z$,
where $V_+$ and $V_-$ are commutative subalgebras.
Recall that $\gt n_\gamma=\gt z$
for generic $\gamma\in X$. In case $\fk=\overline\fk$, one can say immediately that
a generic $(G/N)$-orbit on
$\widetilde X:={\rm Spec}\,\widetilde{\cP}$ has dimension
$(n-d)-(\dim\gt n-1)$. Hence $l(\widetilde X)=l(X)-(\dim\gt n-1)/2$
in case ${\mathcal S}(\gt z)\cap I\ne 0$;
and $l(\widetilde X)=l(X)-(\dim\gt n+1)/2$ otherwise.  
Since the numbers $l(X)$ and $l(\widetilde X)$ do not change under field extensions, 
the same equalities hold over any $\fk$.

If ${\mathcal S}(\gt z)\cap I=0$, then $f_1,\ldots,f_k$ together with
a basis of $V_+\oplus\gt z$ give us a complete
commutative family on $X$. If   
${\mathcal S}(\gt z)\cap I\ne 0$, then we add a basis
of $V_+$ to $\{f_i\}$ and again obtain a complete family on $X$.

$\bullet$ \ If the previous case does not hold, then
either $\gt n$ is a Heisenberg Lie algebra such that
$[\gt n,\gt n]$ {\it is not} a central subalgebra of $\gt g$,
or $\gt n$ contains a commutative characteristic ideal
$\gt h$ such that $\dim\gt h>1$, see Lemma~\ref{nil}.
In both cases, there is a commutative ideal
$\gt h\subset\gt n$ of $\gt g$ such that
either $[\gt g,\gt h]\ne 0$ or $\dim\gt h>1$.
Set $I_0:=I\cap{\mathcal S}(\gt h)$ and let $\gt x_{\gt h}$
be the subvariety of $\gt h^*$ defined by $I_0$.
By definition,
$\gt x_{\gt h}(\overline\fk)$ coincides with the image  of
the natural projection $X(\overline\fk)\to(\gt h\otimes_{\fk}\overline\fk)^*$.  
The connected unipotent subgroup $H=\exp(\gt h)\lhd G$ will play a r\^ole in the proof. 

Set $\mathbb K:=\fk(\gt x_{\gt h})={\rm Quot}\,{\mathcal S}(\gt h)/I_0$.
Consider the  Lie algebra
$\widetilde{\gt g}:=\{\xi\in\gt g\otimes_{\gt h}\mathbb K\mid \{\xi,\gt h\}\subset I_0\otimes\mathbb K\}$
and also set $\widetilde{\cP}:={\cP}^{\gt h}\otimes_{\fk[\gt x_{\gt h}]}\mathbb K$.
By Lemma~\ref{non-closed},  $\widetilde{\cP}$ is a Poisson
quotient of ${\mathcal S}(\widetilde{\gt g})$.

We claim that $\widetilde{\cP}$
contains no zero-divisors. Indeed, suppose  $x,y\in\widetilde{\cP}$ and $xy=0$.
After multiplying $x$ and $y$ by suitable elements
of the field $\fk(\gt x_{\gt h})$, we may assume that
$x,y\in \cP^{\gt h}$.
Since  $\cP$ is a domain, either $x=0$ or $y=0$.

Set $\widetilde X:={\rm Spec}\,\widetilde{\cP}$.
Then $\widetilde X\subset\widetilde{\gt g}^*$ is Poisson
subvariety defined over $\mathbb K$. Let us compute $l(\widetilde X)$. In order to simplify notation, we do
it in case $\fk=\overline\fk$.  (The numbers $l(X)$ and $l(\widetilde X)$ do not change under field extensions.) 

Let $k$ be the dimension of a generic $H$-orbit
on $X$. Note that $k$ is also the dimension
of a generic $G$-orbit in $\gt x_{\gt h}$. Since
$\gt h$ is an algebraic Lie algebra consisting of
nilpotent elements, we have $\fk(X)^H={\rm Quot}\,\fk[X]^H$.
Therefore
generic $H$-orbits on $X$ are separated
by regular $H$-invariants  and ${\rm tr.deg}\,{\cP}^{\gt h}=n-k$. Hence
${\rm tr.deg}\,\widetilde{\cP}=n-k-\dim\gt x_{\gt h}$.

Next, $\mathbb K(\widetilde X)=
\fk(X)^{H}\otimes_{\mathbb K}\mathbb K$.
Recall that $\widetilde{X}$ is a Poisson subvariety of
$\widetilde{\gt g}^*$. In particular, the
Poisson centre ${\eus Z}\mathbb K(\widetilde X)$
 of $\mathbb K(\widetilde X)$ coincides with
$\mathbb K(\widetilde X)^{\widetilde{\gt g}}$.
Because $\gt h$ is commutative, $\gt h\subset \fk[X]^{H}$.
Therefore  the Poisson centre ${\eus Z}\fk(X)^{H}$ is equal to the Poisson centraliser
$$R:=\{f\in\fk(X)\mid \{f,\fk[X]^{H}\}=0\}.$$
Clearly $R$ contains both
$\fk[\gt x_{\gt h}]$ and
${\eus Z}\fk(X)=\fk(X)^{\gt g}$.
For generic $\gamma\in X$ 
we have
$\dim(\gt h_{|\gt g\gamma})=\dim(\gt h\gamma)=k$.
Since all functions in $\fk(X)^{\gt g}$
are constant on $G$-orbits, the subspace of $T^*_\gamma X$ generated
by $d_\gamma\fk[\gt x_{\gt h}]$ and
$d_{\gamma}(\fk(X)^{\gt g})$ has dimension $d+k$. Hence,
${\rm tr.deg}\,R\ge d+k$.
By a simple dimension reason
${\rm tr.deg}\,R= d+k$. Since
${\eus Z}\mathbb K(\widetilde X)=
{\eus Z}\fk(X)^{H}\otimes_{\mathbb F}\mathbb K$,
we get ${\rm tr.deg}\,\mathbb K(\widetilde X)^{\widetilde{\gt g}}=
d+k-\dim\gt x_{\gt h}$. Thus
$$
\begin{array}{r}
l(\widetilde X)=
(\dim_{\mathbb K}\widetilde X+{\rm tr.deg.}\,\mathbb K(\widetilde X)
^{\gt g})/2=
(n-k-\dim\gt x_{\gt h}+d+k-\dim\gt x_{\gt h})/2=\\
 =(n+d)/2-\dim\gt x_{\gt h}=l(X)-\dim\gt x_{\gt h}.
\end{array}
$$

It remains to show that the dimension of $\widetilde{\gt g}$
over $\fk(\gt x_{\gt h})$ is less than $\dim_{\fk}\gt g$.
If this is not the case, then
$\widetilde{\gt g}=\gt g\otimes_{\fk}\mathbb K$
and $\hat{\gt h}=0\ $
(here $\hat{\gt h}$ is the same as in~(\ref{g2})).
From the first equality we get
$[\gt g,\gt h]\subset I_0$, hence $[\gt g,\gt h]=0$;
and by the second one, $\dim\gt h=1$. Together these
conclusions contradict the initial assumptions on $\gt h$.

Applying the inductive hypothesis to $\widetilde X$,
we construct $l(X)-\dim\gt x_{\gt h}$
functions $\tilde f_i\in{\mathcal S}(\widetilde{\gt g})$
such that their restrictions
give us a complete commutative family on $\widetilde X$.
After multiplying them by a suitable
element of $\mathbb K$, we may assume that
$\tilde f_i\in{\mathcal S}(\gt g)$. The remaining
$\dim\gt x_{\gt h}$ functions we get from ${\mathcal S}(\gt h)$.

$\bullet$ \ Suppose now that $\gt g$ is a non-algebraic Lie algebra.
If the nilpotent radical $\gt n\lhd\gt g$ contains a characterisitic ideal $\gt h$ such that
either $\dim\gt h>1$ or $[\gt g,\gt h]\ne 0$, then the above ``commutative" part of the proof 
(decreasing  of $\dim\gt g$) goes without any  alteration.  If $\gt n=0$, then $\gt g$ is reductive and
algebraic. It remains to consider the ``Heisenberg" case.

Choose any decomposition $\gt n=V\oplus\gt z$ and any non-zero $\alpha\in\gt n^*$ with $\alpha(V)=0$.
Then $\tilde{\gt l}:=\{\xi\in \gt g \mid [\xi,V]\subset V\}$ is a subalgebra such that  $\tilde{\gt l}\cap\gt n=\gt z$
and $\tilde{\gt l}+\gt n=\gt g$.  Since $\gt z$ lies in the centre of $\gt g$, the subgroup
$\tilde L=\exp(\tilde{\gt l})\subset G$ can play the r\^ole of $L$ in the proof of Lemma~\ref{heis}.
Here the conclusion  is that  $(\fk[X][1/z])^N\cong \fk[\tilde X][1/z]$, where $\tilde X=X\cap \tilde{\gt l}^*$
is the intersection  in the scheme sense, and  $\tilde{\gt l}^*={\rm Ann}(V)\subset \gt g^*$.
The reduction from $\gt g$ to $\tilde{\gt l}$  still works, because $l(\tilde X)=l(X)-\frac{1}{2}\dim V$.
\end{proof}

\end{document}